\documentclass[11pt]{amsart}
\usepackage{fullpage}
\usepackage[dvips]{graphicx}
\usepackage{amsmath}
\usepackage{psfrag}
\usepackage{amssymb}
\usepackage{epsf}
\usepackage{epsfig}
\newcommand{\prv}[1]{\noindent {\em Proof {#1}} \hspace{0.5cm}}
\newcommand{\fprv}{\hfill\cqfd\myskip }
\def\myskip{\medskip}
\def\cqfd{\nobreak\nopagebreak\rule{0pt}{0pt}\nobreak\hfill\nobreak\rule{.5em}{.5em}}

\newtheorem{theorem}{Theorem}
\newtheorem{lemma}{Lemma}

\newcommand{\GN}{{\mathbb N}}
\newcommand{\GR}{{\mathbb R}}
\newcommand{\GZ}{{\mathbb Z}}
\newcommand{\GC}{{\mathbb C}}
\newcommand{\GK}{{\mathbb K}}
\newcommand{\GQ}{{\mathbb Q}}

\theoremstyle{definition}
\newtheorem{definition}{Definition}

\theoremstyle{remark}
\newtheorem{remark}{Remark}
\numberwithin{equation}{section}
\numberwithin{figure}{section}

\title
{Purely periodic $\beta$-expansions in the Pisot non-unit case} 

\author{Valerie Berth\'e}
\address{LIRMM - 161 rue Ada -34392 Montpellier Cedex 5- France}
\email{berthe@lirmm.fr}

\author{Anne Siegel}
\address{IRISA - Campus de Beaulieu - 35042 Rennes Cedex - France}
\email{Anne.Siegel@irisa.fr}


\begin{document}
\begin{abstract}
It is well known that real numbers with a purely periodic
decimal expansion are the rationals having, when reduced,
a denominator
coprime with $10$.
The aim of this paper is to extend this result  to  beta-expansions with a  Pisot base
beta which is not necessarily a unit: we  characterize
real numbers having  a purely periodic
expansion in such a base; this characterization
is given in terms of an explicit set, called  generalized Rauzy fractal,  
which is shown
to be  a graph-directed self-affine compact subset
of non-zero measure which belongs to 
the   direct product of  Euclidean  and $p$-adic spaces.

\end{abstract}
\maketitle

\bigskip

Let $\beta$ be a Pisot number and $T_{\beta}:x \mapsto \beta x \ (\mbox{mod } 1)$
be the associated $\beta$-transformation.
The aim of this paper is to characterize the real numbers $x$
in $\GQ(\beta) \cap [0,1)$
having a purely periodic $\beta$-expansion. 

It is well known 
that  if $\beta$ is a Pisot number, then the  real numbers that have a
ultimately periodic $\beta$-expansion are the  elements of $\GQ(\beta)$ \cite{bertrand77,schmidt80}.
Thus real numbers $x$ that have a purely
periodic beta-expansion belong to ${\mathbb Q}(\beta)$. 
We present a characterization that 
involves  the conjugates of    the algebraic
 number $x$, and can be compared to Galois' theorem for 
classical continuous fractions.


\begin{theorem}\label{theo:1}
Let $\beta$ be a Pisot number. A real number $ x \in {\mathbb Q}(\beta)\cap [0,1)$ has a purely periodic $beta$-expansion
if and only if $x$ and its conjugates belong to  an explicit subset   in 
 the product of  Euclidean  and $p$-adic spaces (see Figure
\ref{fig_markov} below); this set (denoted $\widetilde{{\mathcal R}_{\beta}}$ and called {\it generalized
Rauzy fractal}) is
a graph-directed self-affine compact subset  in the sense of  \cite{mauldin}
  of non-zero measure;  the primes $p$ that occur are the prime divisors of the norm of $\beta$.
\end{theorem}

The scheme of  the proof  is  based on
a realization of the natural extension of the $\beta$-transformation
 $T_{\beta}$
extended to a geometric space of representation
for the two-sided $\beta$-shift $(X_{\beta},S)$. 
Our results and our proof is  inspired by 
  \cite{Hui2,sano01,SANO-0}  which presents  a similar characterization 
of purely periodic expansions in the case where $\beta$ is a Pisot unit. 

\smallskip

The construction  of the set $\widetilde{{\mathcal R}_{\beta}}$ (introduced in
Theorem \ref{theo:1})
  is inspired
by the geometric representation
as  generalized Rauzy fractals  (also called atomic surfaces)
of substitutive symbolic dynamical systems
in the non-unimod\-ul\-ar case
developed in \cite{sie_ari}. In fact, a substitution $\sigma$ is a non-erasing  morphism 
of the free mono\"{\i}d  ${\mathcal A}^{*}$ and a substitutive dynamical system is a  symbolic
 dynamical system generated by  an infinite sequence which is a fixed point  of a
substitution.
 Furthermore,  if  the $\beta$-expansion of
 $1$ in  base $\beta$   is finite ($\beta$ is said to be a
 simple Parry number) and if its length coincides with the degree
of $\beta$,
then    the set $\widetilde{{\mathcal R}_{\beta}}$ involved in our characterization   is exactly the generalized  Rauzy fractal  that
is associated  in \cite{sie_ari} with the underlying $\beta$-substitution (in the sense of \cite{thurston,Fa}).

Rauzy fractals have been widely studied; for more details, see for instance  \cite{sie_fogg}. There are  mainly 
two methods of construction for   Rauzy fractals.
One  first approach inspired by the seminal paper \cite{rau1},
is  based on formal power series, and is  developed 
in  
 \cite{mes,mes2},  or  in
\cite{can_sie,can_sie1}. 
A second  approach via iterated function systems (IFS)
and generalized substitutions has been developed following ideas from
\cite{itokimura91} in   \cite{arn_ito,sanoarnouxito01,holzam98,SW} 
with special focus on the self-similar properties
of  the Rauzy fractals. 
We combine  here both approaches: we  define the set $\widetilde{{\mathcal R}_{\beta}}$
by  introducing a representation map of the two-sided shift $(X_{\beta}, S)$
 based on  formal power series, and prove that this set  has non-zero Haar measure 
by  cutting  it into pieces that are solutions of an IFS. 
Similar sets  have also been
introduced  and studied  in the framework of $\beta$-numeration
  by 
S. Akiyama in \cite{AKI1,AKI3,AKI6,AKI4},  inspired by
\cite{thurston}.
As an  application of generalized Rauzy fractals, let us mention
that 
they provide Markov partitions for toral 
automorphisms of the torus, as  illustrated in \cite{itooht93,kenyon,Pra99,schmidt00,these_anne}.
Furthermore, there are numerous relations between generalized Rauzy fractals and discrete planes as studied for instance in \cite{ABI,ABS}.

The aim of this paper is twofold. We first want  to characterize 
 real numbers  having a purely periodic  $\beta$-expansion;
we second  try to settle the first steps of a study of  the geometric representation
of $\beta$-shifts in the Pisot non-unit case, generalizing the  results of
\cite{AKI1,AKI3,AKI6,AKI4}, based  on  the formalism
introduced in the  substitutive case in  \cite{sie_ari}.

This paper is organized as follows. We first recall in Section \ref{sec:beta}
the basic elements  needed on $\beta$-expansions. We then 
associate  in Section \ref{sec:rep}
with 
 the two-sided  $\beta$-shift $(X_{\beta},S)$  formal power series in $\GQ[[X]]$;
we  obtain in Section \ref{sec:two}
a representation map  for the two-sided $\beta$-shift
by gathering the set of  finite values 
which can be taken  for any topology (Archimedean or
not) by these  formal power series when specializing them  in $\beta$: in fact, we take
the completion of $\GQ(\beta)$ with respect to
all the absolute values on
$\GQ(\beta)$ which take a value different from
$1$ on $\beta$ (this value is thus smaller than $1$ since $\beta$ is Pisot).
We  are then able to define the Rauzy geometric representation  of the two-sided $\beta$-shift
(Definition \ref{def:rep}). Section \ref{sec:pro} is devoted to 
 the study of the properties   of the set $\widetilde{{\mathcal R}_{\beta}}$.
We then prove Theorem \ref{theo:1} in Section \ref{sec:car}.

\section{$\beta$-numeration}\label{sec:beta}

Let $\beta >1$ be a real number. {\it In all that follows, $\beta$ is assumed to be a Pisot number}.
 The Renyi $\beta$-expansion of a real number $x \in [0,1)$
is defined as the sequence
$(x_i)_{i\geq 1}$ with values in ${\mathcal A}_{\beta}:= \{0, 1, \dots, 
[ \beta ]\}
$ produced by
the $\beta$-transformation 
$T_{\beta}: \ x \mapsto \beta x \ (\mbox{mod } 1)\ $
as follows
$$\forall i \geq 1, \ u_i=\lfloor \beta T_{\beta}^{i-1}(x)\rfloor,
\mbox{ and thus }
x =\sum_{i\geq 1} u_i \beta^{-i}.$$
Let $d_{\beta}(1)=(t_i)_{i\geq 1}$ denote the $\beta$-expansion of $1$.
Numbers $\beta$ such that $d_{\beta}(1)$ is ultimately periodic
are called {\em Parry numbers} 
and those such that
$d_{\beta}(1)$ is finite        are called {\em simple Parry numbers}. 
Since $\beta$ is assumed to be Pisot, then 
$\beta$ is either a  Parry number or a simple Parry number \cite{BM86}.
 Let $d^*_{\beta}(1)=d_{\beta}(1)$, if $d_{\beta}(1)$
is infinite, and $d^*_{\beta}(1)=(t_1 \dots t_{n-1} (t_n -1))^{\infty}$,
if $d_{\beta}(1)=t_1\dots t_{n-1}t_n$ is finite ($t_n \neq 0$). 
The set of $\beta$-expansions of real numbers in $[0,1)$
is exactly the set  of sequences $(u_i)_{i\geq 1}$ in
 ${\mathcal A}_{\beta}^{\GN}$
such that
\begin{equation}\label{eq:real}
\forall k \geq 1, \ (u_i)_{i\geq k} <_{\mbox{lex}} d^*_{\beta}(1).
\end{equation}

For more details on the $\beta$-numeration, 
see for instance \cite{Loth,frou_temuco}.

\subsubsection*{The (two-sided symbolic) $\beta$-shift}
Let  $(X_{\beta}, S)$ 
denote the two-sided symbolic dynamical system
associated with $\beta$, where the shift map $S$ maps
the sequence $(y_i)_{i\in \GZ}$ onto $(y_{i+1})_{i\in \GZ}$.
The set $X_{\beta}$  is defined as  the set of two-sided sequences $(y_i)_{i\in \GZ}$ in ${\mathcal A}_{\beta}^{\GZ}$
such that each left truncated sequence is less
than or equal 
to $d^*_{\beta}(1)$, that is,
$\forall k \in \GZ, \ (y_i)_{i\geq k} \leq _{\mbox{lex}} d^*_{\beta}(1).$

We will use the following notation for 
the elements of $X_{\beta}$:
if $y=(y_i)_{i\in \GZ} \in X_{\beta}$, define
$u=(u_i)_{i\geq 1}=(y_i)_{i\geq 1}$
and $w=(w_i)_{i\geq 0}=(y_{-i})_{i\geq 0}$.
One thus gets a two-sided sequence of the form
$$\dots w_3w_2w_1w_0u_1u_2u_3\dots$$
and  write it as $y=((w_i)_{i\geq 0},(u_i)_{i\geq 1})=(w,u).$
In other words, we will use the letters
$(w_i)$ for denoting  the ``past'' and $(u_i)$ for the ``future''
of the element $y=(w,u)$
of the two-sided shift $X_{\beta}$.

\subsubsection*{One-sided $\beta$-shifts}
We denote by $X_{\beta}^r$  the set of one-sided sequences
$u=(u_i)_{i\geq 1}$ such that
there exists $w=(w_i)_{i\geq 0}
$ with $(w,u) \in X_{\beta}.$ This set is  called the
{\em right one-sided $\beta$-shift}.  
        It coincides with the usual one-sided $\beta$-shift and is equal to the set of sequences 
$(u_i)_{i\geq 1}$ which satisfy
\begin{equation}\label{eq:real2}
\forall k \geq 1, \ (u_i)_{i\geq k} \leq _{\mbox{lex}} d^*_{\beta}(1).
\end{equation}
One similarly defines $X_{\beta}^l$  as the set of one-sided sequences
$w=(w_i)_{i\geq 0}$ such that
there exists $u=(u_i)_{i\geq 1}
$ with $(w,u) \in X_{\beta}$.  We call it the {\em left one-sided
 $\beta$-shift}.

\subsubsection*{Sofic shift}
\begin{figure}[ht]
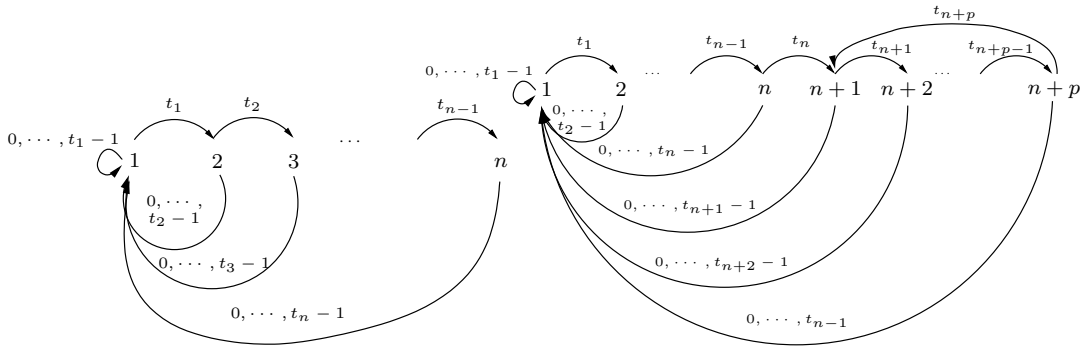

\begin{center}
\input{graph1.pstex_t}\input{graph2.pstex_t}
\caption{The  automata ${\mathcal M}_{\beta}$ for $\beta$ simple Parry number ($d_{\beta}(1)=t_1\dots t_{n-1}t_n$)
 and for $\beta$ non-simple Parry number ($d_{\beta}(1)=t_1\dots t_n(t_{n+1}\dots t_{n+p})^{\infty}$).}\label{fig_graphe}\end{center}
\end{figure}

Since
$\beta$ is  a  Parry number (simple or not), then $(X_{\beta}, S)$ is
sofic \cite{BM86}.
We denote by $F(X_{\beta})$ the set of finite factors
of the sequences in $X_{\beta}$;
the minimal automaton ${\mathcal M}_{\beta}$ recognizing the set of factors
of  $F(D_{\beta})$
can easily be constructed (see Figure \ref{fig_graphe}).  The number of states
$d$ of this automaton is equal to   the  length of the period $n$ of $d^*_{\beta}(1)$
if $\beta$ is a simple Parry number with $d_{\beta}(1)=t_1\dots t_{n-1}t_n$, $t_n \neq 0$,
and to the sum of its preperiod $n$ plus its period $p$,
if $\beta$ is a non-simple Parry number with
 $d_{\beta}(1)=t_1\dots t_n(t_{n+1}\dots t_{n+p})^{\infty}$ ($t_n \neq t_{n+p}$, $t_{n+1}\cdots t_{n+p} \neq 0^p$).

\subsubsection*{$\beta$-substitutions}
Let us recall that a  \emph{substitution} $\sigma$ is a morphism of the free mono\"{\i}d  ${\mathcal A}^{*}$, 
such that the image of each letter
of  $ {\mathcal A} $ is non-empty. As introduced for instance in \cite{thurston} and in \cite{Fa},
one can associate in a natural way
a substitution  $\sigma_{\beta}$ (called $\beta$-substitution)
  with $(X_{\beta}, S)$ 
over the alphabet $\{1,\cdots,d\}$, where $d$ denotes the number of states
of the automaton ${\mathcal M}_{\beta}$: $j$ is the $k$-th letter occuring in $\sigma_{\beta}(i)$ (that is,
$\sigma_{\beta}(i)=pj s$, where $p,s \in \{1,\cdots,d\} ^*$
and $|p|=k-1$) if and only if there is an arrow in ${\mathcal M}_{\beta}$ from the state
$i$ to the state $j$ labeled by $k-1$. One easily checks that this definition is consistent.

An explicit formula for $\sigma_{\beta }$ can be computed by considering
the  two different cases, $\beta$ simple and $\beta$ non-simple Parry number. 
\begin{itemize}
\item Assume $d_{\beta}(1)=t_1\dots t_{n-1}t_n$ is finite, with $t_n\neq 0$.
Thus $d^*_{\beta}(1)=(t_1\dots t_{n-1}(t_n-1))^{\infty}$.
One defines $\sigma_{\beta}$ over the alphabet
$\{1,2, \dots,n\}$  as shown in  (\ref{formuleSub}). 

\item Assume $d_{\beta}(1)$ is infinite. Then it cannot be purely periodic (according to
Remark 7.2.5 \cite{Loth}). Hence
$d_{\beta}(1)=d^*_{\beta}(1)=t_1\dots t_n(t_{n+1}\dots t_{n+p})^{\infty}$,
with $n\geq 1$, $t_n\neq t_{n+p}$ and $t_{n+1}\dots t_{n+p} \neq 0^p$.
One defines $\sigma_{\beta}$ over the alphabet
$\{1,2, \dots,n+p\}$  as shown in  (\ref{formuleSub}). 
\end{itemize}

\begin{equation}\label{formuleSub}
\begin{array}{ccc}
\sigma_{\beta}:\left\{\begin{array}{ll}
1 &\mapsto 1^{t_1}2 \\
2&\mapsto 1^{t_2}3\\
\vdots&\vdots \\
n-1 &\mapsto 1^{t_{n-1}}n\\
n &\mapsto 1^{t_n}.
\end{array}\right.   
& & 
\sigma_{\beta}:\left\{\begin{array}{ll}
1 &\mapsto 1^{t_1}2 \\
2&\mapsto 1^{t_2}3\\
\vdots&\vdots \\
n+p-1&\mapsto 1^{t_{n+p-1}}(n+p)\\
n+p &\mapsto 1^{t_{n+p}}(n+1).
\end{array}\right. 
\\
& \hspace{1cm}& \\
\mbox{Substitution associated with } & & \mbox{Substitution associated with } \\
\mbox{a simple Parry number} & &  \mbox{a non-simple Parry number}
\end{array}
 \end{equation}

 We will use $\beta$-substitutions in Section \ref{sec:pro}
in order to describe properties of the set $\widetilde{{\mathcal R }_{\beta}}$. 
Notice that the automaton  ${\mathcal M }_{\beta}$ is exactly the prefix-suffix automaton 
 of the substitution 
$\sigma_{\beta}$  considered  for instance in \cite{can_sie,can_sie1},
 after reversing all the edges and replacing the labels in the prefix-suffix automaton
by the lengths of the prefixes.

 The {\it incidence matrix} of the substitution $\sigma_{\beta }$ is defined as
the $d\times d$ matrix whose entry of index $(i,j)$ counts the number
of occurrences of the letter $i$ in $\sigma (j)$.  As a consequence of the definition, the incidence
matrix of $\sigma_{\beta }$ coincides with the transpose of the adjacency
matrix of the automaton ${\mathcal M }_{\beta}$. By eigenvalue of a substitution
$\sigma$, we  mean in all that follows an eigenvalue of the characteristic polynomial
of the  incidence matrix of $\sigma$. A substitution is said to be of {\em Pisot type} if all
its eigenvalues except its largest one  which is assumed to be simple,  are  non-zero and of  modulus  smaller than 1.

\section{Representation of $\beta$-shifts}\label{sec:rep}

The right  one-sided  shift $X_{\beta}   ^r$
admits as a natural geometric representation
the interval $[0,1]$; namely, one associates with a sequence
$(u_i)_{i\geq 1} \in X_{\beta}^r$ its real value
$\sum_{i\geq 1} u_i \beta^{-i}$.  We even have a measure-theoretical 
 isomorphism
between $X_{\beta}^r$ endowed with the shift, and 
$[0,1]$ endowed with the map $T_{\beta}$.
We want now to give a  similar geometric interpretation
of the set $X_{\beta}$; for that
purpose,  we first give a geometric representation of  $X_{\beta}^l$
as a  generalized Rauzy fractal.

\subsection{Representation of the left one-sided shift $X_{\beta}^l$}\label{sec:left}
The aim of this section is to introduce
first a formal power series, called formal  representation of $X_{\beta}^l$ 
and second a geometrical representation as an explicit compact set of
in 
 the product of Euclidean and $p$-adic spaces 
 following \cite{sie_ari}; the primes which appear as 
$p$-adic spaces
here 
will be  the prime factors of the norm of $\beta$.


\subsubsection*{Formal representation of the symbolic dynamical
 system $X_{\beta}^l$} 

\begin{definition}
The {\em formal representation} of $X_{\beta}^l$ is
denoted $${\varphi}_X: X_{\beta} \to {\mathbb Q}[[X]]$$
where ${\mathbb Q}[[X]]$  is the ring
of formal power series with coefficients in $\GQ$, and defined by:
$$
\mbox{for all } (w_i)_{i\geq 0} \in X_{\beta}^l, \ \  
 {\varphi}_X(w_i) = \sum_{i \geq 0} w_i X^{i}\in \GQ[[X]].
$$

\end{definition}

\subsubsection*{Topologies over $\GQ(\beta)$}

We now  want to specialize  these formal power series by giving to the indeterminate 
$X$ the value $\beta$, and associating with them
values by making them converge.
 We thus  want  to  find a topological framework in which all the series
$\sum _{i \geq 0} w_i \beta ^i$  for $(w_i)_{i\geq 0} \in X_{\beta}^l$,
would converge; in fact, this boils down to find  
 all the Archeme\-dean and non-Archemedean metrizable 
topologies on $\GQ(\beta)$  
for which  these series  converge in a suitable completion;   they are of two types. 
\begin{itemize}
\item  Suppose that the topology (with absolute value $|\cdot|$)
is Archemedean:
its restriction to ${\mathbb Q}$ corresponds to the  usual
absolute value on ${\mathbb Q}$ and there exists a $\GQ$-isomorphism $\tau_i$ such that 
$|x|=|\tau_i(x)|_{\GC}$, for $x \in  \GQ(\beta)$.
The series ${\varphi}_X$ specialized in $\beta$
converge in $\GC$ if and only if
$\tau_i$ is associated with  a conjugate $\beta_i$
of modulus strictly smaller than one.
\item  Assume  that the topology is non-Archimedean: 
there exists a prime ideal  ${\mathcal I}$  of the integer ring
${\mathcal O}_{\GQ(\beta)}$
of ${\mathbb Q}(\beta)$  for which the topology coincides with the
the ${\mathcal I}$-adic topology;
let  $p$  be the prime number defined by
${\mathcal I}\cap \GZ=p\GZ$;
the restriction 
of the topology to ${\mathbb Q}$ is the $p$-adic topology. 
 The series ${\varphi}_X$ specialized in $\beta$
take  finite values in the completion ${\mathbb K}_{{\mathcal I}}$ 
 of ${\mathbb Q}(\beta)$
for the ${\mathcal I}$-adic topology 
if and only if $\beta \in {\mathcal I}$, i.e., 
$|\beta|_{\mathcal I}<1$.
\end{itemize}

\subsubsection*{Representation space ${\mathbb K}_{\beta}$
of $X_{\beta}^l$}
We assume now that $\beta$ be a Pisot number of degree $d$, say.
Let  $\beta_2$, $\dots$, $\beta_r$ be the 
real conjugates of $\beta$ (they all have modulus strictly smaller than $1$,
since $\beta$ is Pisot),
and 
let ${\beta}_{r+1}$, $\overline{\beta_{r+1}}$, $\dots$, $\beta_{r+s}$,
 $\overline{\beta_{r+s}}$ be
its complex conjugates.
For $2 \leq j \leq  r$, let ${\mathbb K}_{\beta_j}$ be equal
to ${\mathbb R}$, and for $r+1\leq j \leq
r+s$, let ${\mathbb K}_{\beta_j}$ be equal
to  ${\mathbb C}$, ${\mathbb R}$
 and ${\mathbb C}$  being endowed with
the usual topology.  For $i=1$ to $d$, let $\tau_i$ be a  $\GQ$-automorphism 
of $\GK=\GQ(\beta_1,\dots,\beta_d)$
which sends $\beta$ on its
algebraic conjugate
$\beta_i$.  For  a given 
$i$ and
for every element $Q(\beta)$ of $\GQ(\beta)$, then $\tau_i(Q(\beta))=Q(\beta_i)$.

We first  gather the complex representations  by omitting the ones which are conjugate
in the complex case.
This representation contains all
the possible Archemedean
values for
${\varphi}_X$.
It takes values in 
$${\mathbb K}_{\infty}={\mathbb K}_{\beta_2}\times \dots\times{\mathbb K }_{\beta_{r+s}}.$$
Let ${\mathcal I}_1$, $\dots$,  ${\mathcal I}_{\nu}$
be  the prime ideals in  the integer ring ${\mathcal O}_{\GQ(\beta)}$
of ${\mathbb Q}(\beta)$ that contain $\beta$, that is,
 \begin{equation}\label{decompo}
\beta {\mathcal O}_{\GQ(\beta)}
=\prod_{i=1}^\nu {\mathcal I_i}^{n_i}.
\end{equation}

Recall that ${\mathbb K}_{{\mathcal I}}$ denotes the completion of ${\mathbb Q}(\beta)$ for
 the ${\mathcal I}$-adic topology.
We then gather the representations in the completions of
$\GQ(\beta)$ for the non-Archemedean  topologies. Hence one defines the  {\em representation space of $X_{\beta}$}  as the
 direct product  ${\mathbb K}_{\beta}$ of
all these fields:
$${\mathbb K}_{\beta}= {\mathbb K}_{{\beta}_2} \times
\dots {\mathbb K}_{{\beta}_{r+s}} \times {\mathbb K}_{{\mathcal I}_1} \times \dots
{\mathbb K}_{{\mathcal I}_{\nu}}\simeq {\mathbb R}^{r-1} \times 
{\mathbb C}^{s}\times {\mathbb K }_{{\mathcal I_1}}
\times \dots \times {\mathbb K }_{{\mathcal I_\nu}}
.$$
The field ${\mathbb K}_{\mathcal I}$
is a  finite extension  of the $p_{\mathcal I}$-adic field
$\GQ_{p_{\mathcal I}}$ where ${\mathcal I}
\cap \GZ=p_{\mathcal I}\GZ.$
For a given prime $p$,
the fields ${\mathbb K}_{\mathcal I}$ that are
$p$-adic fields are the ones for which ${\mathcal I}$
contains simultaneously $p$ and $ \beta$.
Furthermore,  the prime numbers 
$p$ for which there exists a prime ideal
of ${\mathcal O}_{\GQ(\beta)}$ which
contains simultaneously $p$ and $ \beta$
are exactly the prime divisors
of the constant term of the minimal polynomial of $\beta$ (see Lemma 4.2  \cite{sie_ari}).
In particular, ${\mathbb K}_{\beta}$ is  a Euclidean space 
if and only if $\beta$ is a unit. Endowed with the product of
the topologies of each of its elements, ${\mathbb K}_{\beta}$
is a metric abelian group.

 The {\em canonical embedding} of ${\mathbb Q}(\beta)$ into
${\mathbb K}_{\beta}$ is defined by the following  morphism:
\begin{equation}\label{eq:delta}
\delta_{\beta}: P(\beta) \in {\mathbb Q}(\beta) \mapsto
(\underbrace{P(\beta_2)}_{\in{\mathbb K}_{{\beta}_2}} , \dots,
 \underbrace{P(\beta_{r+s})}_{\in{\mathbb K}_{{\beta}_{r+s}}},
\underbrace{P(\beta)}_{\in{\mathbb K}_{{\mathcal I}_1}}, 
\dots, \underbrace{P(\beta)}_{\in{\mathbb K}_{{\mathcal I}_{\nu}}})
 \in {\mathbb K}_{\beta}.
\end{equation} 

Since the  topology on ${\mathbb K}_{\beta}$ has been chosen
so that the formal power series 
$
\lim_{n \to \infty} \delta_{\beta}(\sum_{i=0}^n w_i \beta^i)$ $
=\sum_{i\geq 0} w_i \delta_{\beta}(\beta)^i$
are  convergent in ${\mathbb K}_{\beta}$  for  every $(w_i)_{i\geq 0} \in X_{\beta}^l$, one 
defines the following, where the notation 
$\delta_{\beta}(\sum_{i\geq 0} w_i \beta^i)$
stands for $\sum_{i\geq 0} w_i \delta_{\beta}(\beta)^i$.
\begin{definition}\label{def:map}
The {\em representation map of $X_{\beta}^l$}, called
{\em one-sided representation map},  is
defined by
$$ \varphi_{\beta}:\  X_{\beta}^l
\rightarrow {\mathbb K}_{\beta}, \ 
(w_i)_{i \geq 0}   \mapsto 
\delta_{\beta}(\sum_{i\geq 0} w_i \beta^i)
.$$
We set ${\mathcal R}_{\beta}:=\varphi_{\beta}(X_{\beta}^l)$ 
and call it the {\em generalized
Rauzy  fractal} or {\em 
geometric representation}
 of  the left one-sided shift $X_{\beta}^l$.
\end{definition}

\subsection{Examples}\label{sec:ex}
\subsubsection*{The golden ratio}
Let $\beta= (1+\sqrt{5}) / 2$ be the golden ratio, that is, the largest root
of $X^2-X-1$. One has $d_{\beta}(1)=11$ ($\beta$ is a  simple Parry number)
and $d^*_{\beta}(1)=(10)^{\infty}$. Hence $X_{\beta}$ is
the set of sequences in $\{0,1\}^{\mathbb Z}$ in which
there are no two consecutive 1's.   Furthermore,  the associated $\beta$-substitution is
the Fibonacci substitution: $\sigma_{\beta}: \ 1 \mapsto 12,
  \ 2 \mapsto 1.$ 
One has $\GK_{\beta}=\GR$; the canonical
embedding $\delta_{\beta}$ is reduced to the map $\tau_{(1-\sqrt{5}) / 2}$
(that is, the $\GQ$-automorphism
of $\GQ(\beta)$
 which  maps
$\beta$ on its conjugate),
and 
$\delta_{\beta}(\GQ(\beta))=\GQ(\beta)$.
The set ${\mathcal R}_{\beta}$
is an interval.

\subsubsection*{The Tribonacci number}
Let $\beta$ be the {\em Tribonacci number}, that is, 
the Pisot root of the  polynomial $X^3 - X^2 - X - 1$.
One has $d_{\beta}(1)=111$ ($\beta$ is a  simple Parry number)
and $d^*_{\beta}(1)=(110)^{\infty}$. Hence $X_{\beta}$ is
the set of sequences in $\{0,1\}^{\mathbb Z}$ in which
there are no three consecutive 1's.  Furthermore, $\sigma_{\beta}: \ 1 \mapsto 12,
  \ 2 \mapsto 13, \  3 \mapsto 1.$ 
One has $\GK_{\beta}=\GC$;
the canonical embedding is reduced
to the $\GQ$-isomorphism $\tau_{\alpha}$ which  maps $\beta$
on $\alpha$, where 
$\alpha$ is one of the complex roots  of $X^3 - X^2 - X - 1$.
The set ${\mathcal R}_{\beta}$ which satisfies
$${\mathcal R}_{\beta}=\{\sum_{i\geq 0} w_i\alpha^i;\  \forall i, \ w_i\in\{0,1\}, \ w_iw_{i+1}w_{i+2} \neq 0\}$$
 is a
compact subset of ${\mathbb C}$ called the {\em Rauzy fractal}.
This set was introduced in \cite{rau1}, see also
\cite{itokimura91,mes,mes2}.
It is shown in Fig. \ref{fig_fractal} with a division into
three pieces indicated by different shades.
They correspond to the sequences $(w_i)_{i \geq 0}$
such that either $w_0=0$, or $w_0w_1=10$, or  $w_0w_1=11$. 
There are as many pieces as  the length of
$d_{\beta}(1)$, which is also  equal here to  the degree of
$\beta$. We will come back  to  the interest of this  division of the
Rauzy fractals into smaller pieces in Section \ref{sec:pro}.

\subsubsection*{The smallest Pisot number}
Let $\beta$ be the Pisot root of $X^3-X-1$.
One has  $d_{\beta}(1)=10001$ ($\beta$ is a  simple Parry number)
and $d^*_{\beta}(1)=(10000)^{\infty};$ $\sigma_{\beta}: \ 1 \mapsto 12$,
$2 \mapsto 3, \ 3 \mapsto 4, \ 4 \mapsto 5, \ 5 \mapsto 1$; the characteristic polynomial
of its incidence matrix is
$(X^3-X-1)(X^2-X+1)$, hence $\sigma_{\beta}$ is not a substitution of Pisot type.
A  self-similar tiling generated by it
has been studied in details in \cite{AKI3};
some connected  surprising tilings have also  been introduced
in \cite{ItoEi}.
One has $\GK_{\beta}=\GC$; the canonical embedding is also reduced
to the $\GQ$-isomorphism $\tau_{\alpha}$ which  maps $\beta$
on $\alpha$, where 
$\alpha$ is one of the complex roots of $X^3-X-1$.
 The set ${\mathcal R}_{\beta}$
is shown in Fig. \ref{fig_fractal} with a division into
five pieces  corresponding respectively  to the sequences $(w_i)_{i \geq 0}$
such that either $w_0w_1w_2w_3=0000$, or $w_0=1$,  $w_0w_1=01$,
$w_0w_1w_2=001$, or $w_0w_1w_2w_3=0001$. 
 The number
of different pieces is equal to the
length  of $d_{\beta}(1)$;  there are here
 $5$ pieces  whereas the degree of $\beta$ is $3$.

\subsubsection*{A non-unit example}
\begin{figure}[ht]
\begin{center}
\begin{minipage}{3.5cm}
\includegraphics[height=3.5cm]{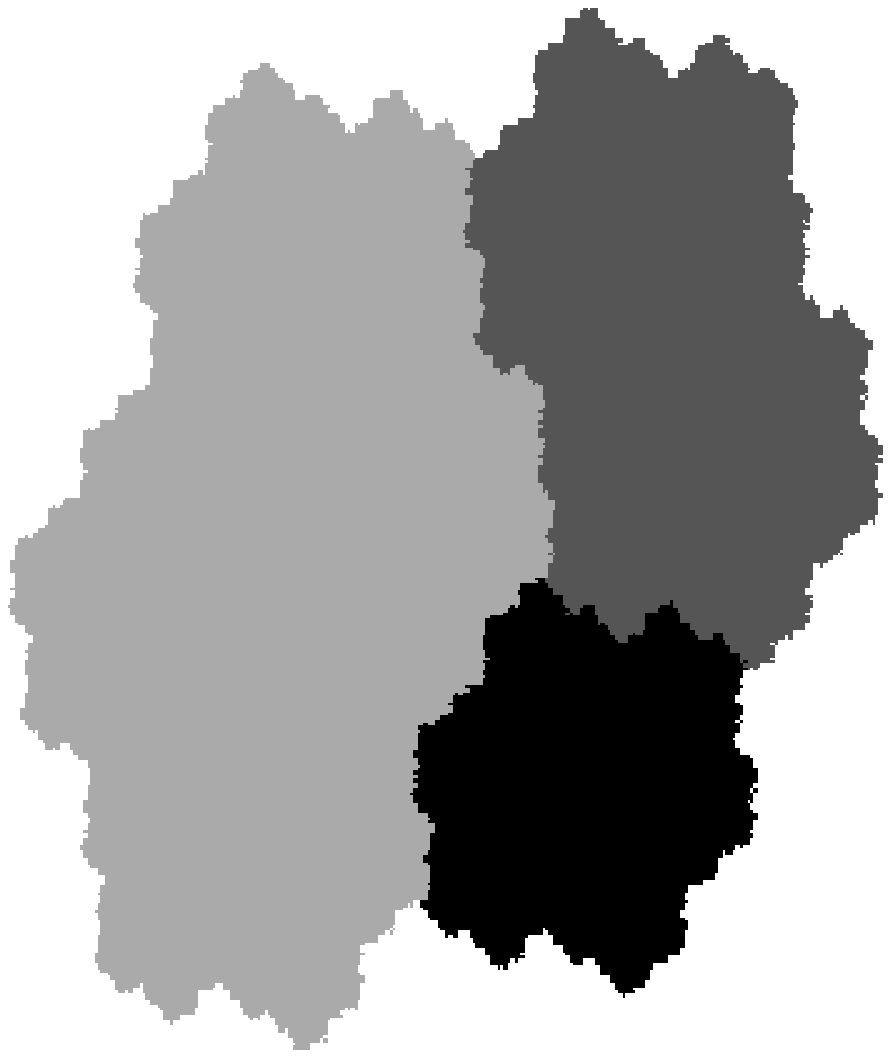}
\end{minipage}
\begin{minipage}{5cm}
\includegraphics[height=3.5cm]{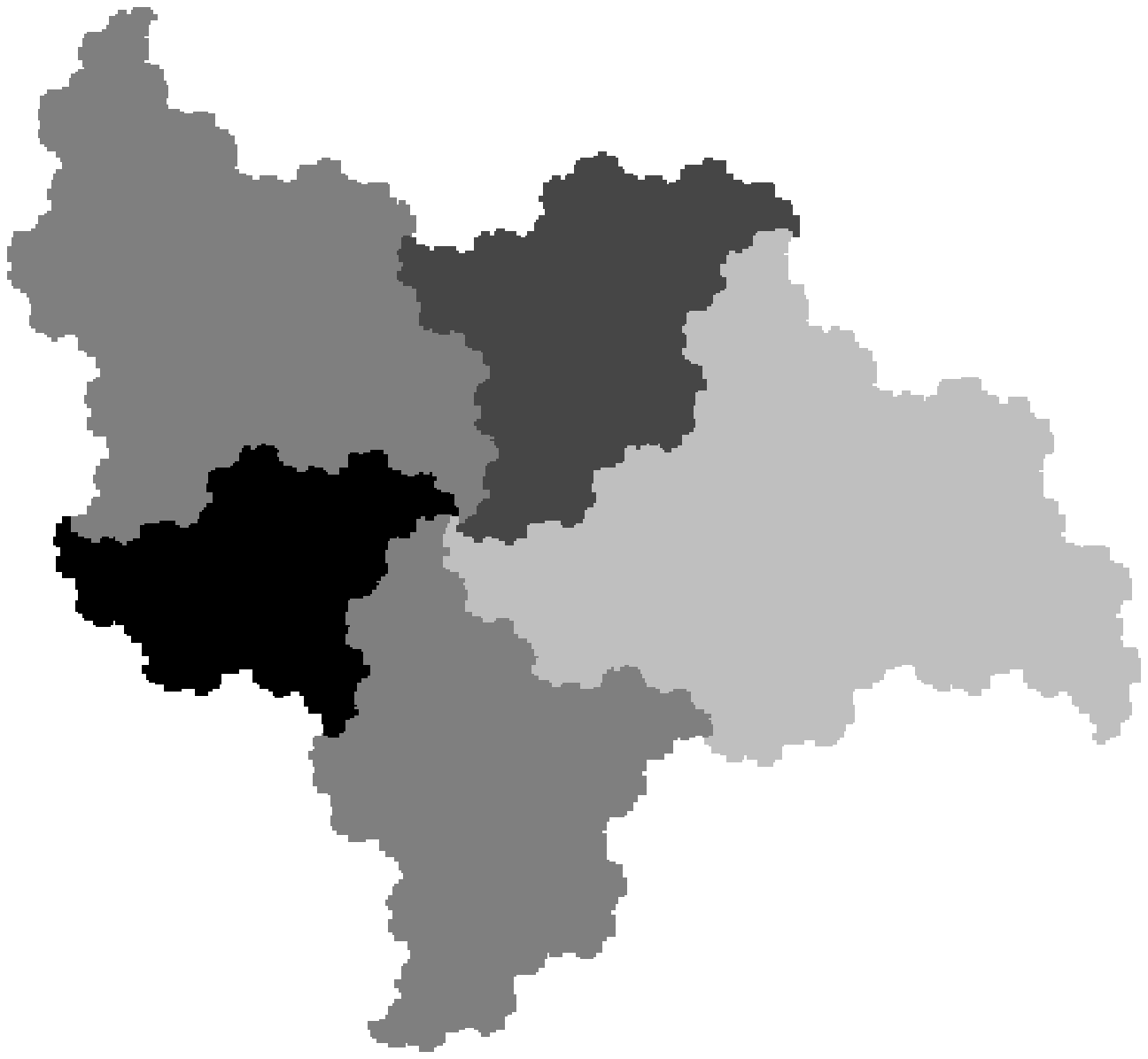}
\end{minipage}
\begin{minipage}{3.5cm}
\includegraphics[height=3.5cm]{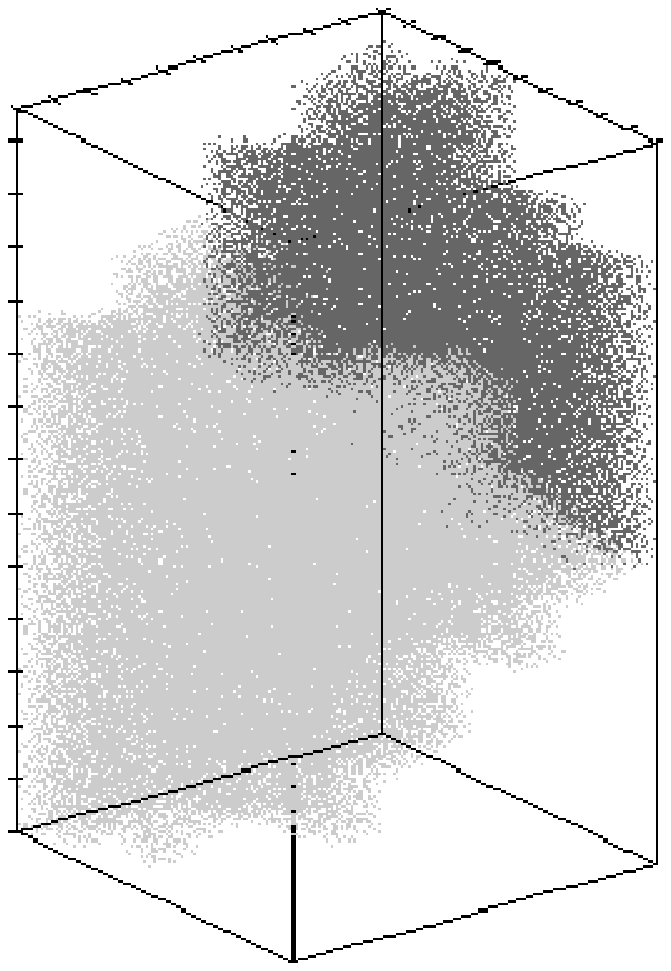}
\end{minipage}
\caption{The Rauzy geometric one-sided representation  for the Tri\-bo\-nacci-shift, the smallest Pisot number-shift and the $(2+\sqrt2$)-shift.}
\end{center}\label{fig_fractal}
\end{figure}

Let $\beta=2 + \sqrt2$  be the dominant root of
the   polynomial  $X^2-4X+2$. The other root is $ 2 - \sqrt 2$.
One has $d_{\beta}(1)=d^*_{\beta}(1)=31^{\infty}$; $\beta$ is not a  simple Parry number; $\sigma_{\beta}: \ 1 \mapsto 1112$,
$2 \mapsto 12$. The ideal
$2\GZ$ is ramified in ${\mathcal O}_{\GQ(\sqrt2)}$, that is,
$2\GZ={\mathcal I}^2$.
Hence there exists only one ideal which contains
$\sqrt 2$; its index of ramification is $2$;
the degree of the extension  $\GK_{\mathcal I}$ 
over $\GQ_2$ has degree $2$.
Hence the  geometric one-sided  representation 
${\mathcal R}_{\beta}$ is a subset
of $\GR \times\GQ_2\times \GQ_2$, shown in Fig. \ref{fig_fractal}. 
The division of the Rauzy fractal 
cannot be expressed as in the previous examples as    finite conditions
 on the prefixes of the sequences $(w_i)_{i\geq 0}$ (for more details,
see Section \ref{sec:ex2}).

\subsection{Representation of the two-sided shift
$(X_{\beta},S)$}\label{sec:two}
 
 \subsubsection*{Representation space $\widetilde{\GK_{\beta}}$ of the two-sided $\beta$-shift $X_{\beta}$
and  representation map $\widetilde{\varphi_{\beta}}$} 
We define now 
the {\em representation map} $\widetilde{\varphi_{\beta}}$
of $X_{\beta}$: it takes its values 
 in 
$\GK_{\beta}\times \GR$  and maps   a point
 $(w,u)=((w_i)_{i\geq 0},(u_i)_{i\geq 1})$
in $X_{\beta}$  to the point  obtained 
by  introducing  on the last
 coordinate (in $\GR$)
the real number whose $\beta$-expansion  is given by $(u_i)_{i\geq 1}$,
and  by
gathering on the first coordinate (in $\GK_{\beta}$)
the set of  finite values  which can be taken by  
the formal power series ${\varphi}_X(w)$  specialized in $\beta$ for all the topologies 
that exist on ${\mathbb Q}(\beta)$. 
Let us define $\widetilde{{\mathbb K}_{\beta}}$ as 
${\mathbb K}_{\beta}\times \GR$.

\begin{definition}\label{def:rep}
The {\em representation map $\widetilde{\varphi_{\beta}}: X_{\beta} \rightarrow  \widetilde{{\mathbb K}_{\beta}}$ of $X_{\beta}$}, called
{\em two-sided representation map},  is
defined for all  $(w,u)=((w_i)_{i\geq 0}, (u_i)_{i\geq 1}) \in X_{\beta}$ by:
$$
\widetilde{\varphi_{\beta}}(w,u)  =  (- \varphi_{\beta}(w_i),\sum_{i\geq 1} u_i \beta^{-i})
=(-\delta_{\beta}(\sum_{i\geq 0}w_i \beta^{i}),
\sum_{i\geq 1}  u_i{{\beta}}^{-i}) 
  \in{\mathbb K}_{\beta}  \times \GR.
$$
The set $\widetilde{{\mathcal R}_{\beta}}:=\widetilde{\varphi_{\beta}}(X_{\beta})$ is called
the {\em Rauzy  fractal}
or {\em geometric representation }
 of the two-sided $\beta$-shift.  This set  is easily seen to be bounded and hence compact.
\end{definition}

We will see below (see (\ref{eq:prop1}) and (\ref{eq:minus}))
 the interest of introducing the sign minus before $\delta_{\beta}$ in the definition of the map
$\widetilde{\varphi_{\beta}}$.

\subsubsection*{Extension $\widetilde{T}_{\beta}$ of the $\beta$-transformation.}
One can extend in a natural way the definition
of $T_{\beta}$ to the product of the representation space
${\mathbb K}_{\beta}$  by $\GR$ as follows, in order to obtain
a realization of its natural extension.

\begin{definition}\label{defi:h}
Let  $h_{\beta}: \ {\mathbb K}_{\beta}
\rightarrow {\mathbb K}_{\beta} $
stands  for the multiplication  map in ${\mathbb K}_{\beta}$
by the diagonal matrix whose diagonal  coefficients are given  by  $\delta_{\beta}(\beta)$.
 One thus defines
$$\widetilde{T_{\beta}}: \, {\mathbb K}_{\beta}
\times \GR
\rightarrow {\mathbb K}_{\beta}
\times \GR, \ (a,b)
  \mapsto ( h_{\beta}({ a}) - [\beta b] \delta_{\beta}(1),\beta b - [\beta b]).$$
\end{definition}

In particular, the following commu\-ta\-tion re\-la\-tions hold,
whe\-re $Id$ denotes the identity map
over $\GR$;  recall that $S$ denotes the shift over $X_{\beta}$: 
\begin{equation}\label{eq:prop2}
\widetilde{T_{\beta}} \circ 
(\delta_{\beta},Id)=(\delta_{\beta},Id) \circ T_{\beta}
\mbox{ over } {\mathbb Q}(\beta).
\end{equation}
\begin{equation} \label{eq:prop1}
\forall \, (w,u) \in X_{\beta}, \quad  \mbox{ if } u\not= d^*_{\beta}(1), \quad \widetilde{\varphi_{\beta}} \circ S (w,u)= \widetilde{T_{\beta}}  \circ
 \widetilde{\varphi_{\beta}}(w,u).
\end{equation}

The proof of (\ref{eq:prop2}) is immediate.
Let us prove  (\ref{eq:prop1}). Take an element   $(w,u)=((w_i)_{i\geq 0},(u_i)_{i\geq 1}) \in {X_{\beta}}$
with  $u\not= d^*_{\beta}(1)$.
One has  $\widetilde{\varphi_{\beta}}(w,u)  =(-\delta_{\beta}(\sum_{i\geq 0}w_i \beta^{i}),
\sum_{i\geq 1}  u_i{{\beta}}^{-i})$. 
Since $u$ satisfies (\ref{eq:real}), then  the integer part of $\beta (\sum_{i\geq 1}  u_i{{\beta}}^{-i})$
is equal to $u_1$,
hence $T_{\beta}( \sum_{i\geq 1}  u_i{{\beta}}^{-i})= \sum_{i\geq 2}  u_i{{\beta}}^{-i}$,
and the assertion follows. 

\subsubsection*{Examples} 
\begin{figure}[ht]
\begin{center}
\begin{minipage}{3.1cm}
\includegraphics[height=3cm]{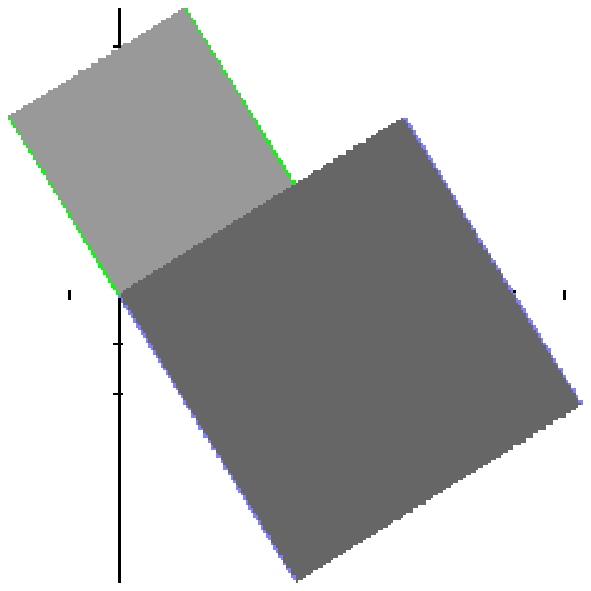}
\end{minipage}
\begin{minipage}{4.1cm}
\includegraphics[height=4cm]{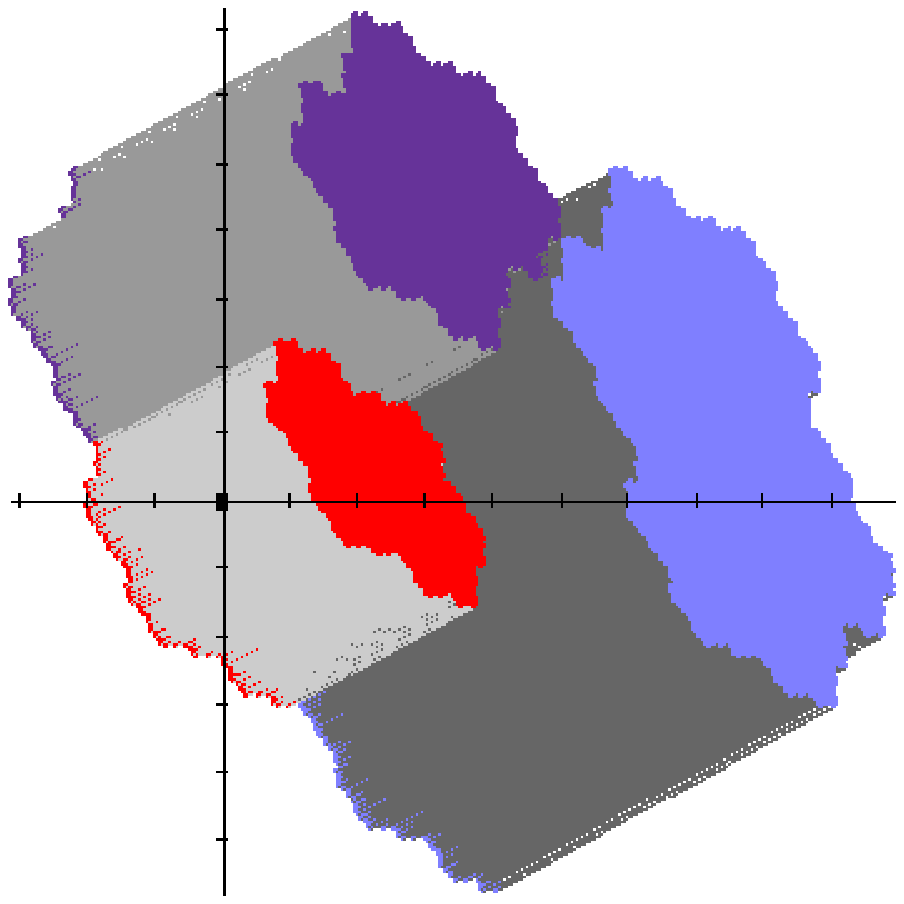}
\end{minipage}
\begin{minipage}{5.1cm}
\includegraphics[height=4cm]{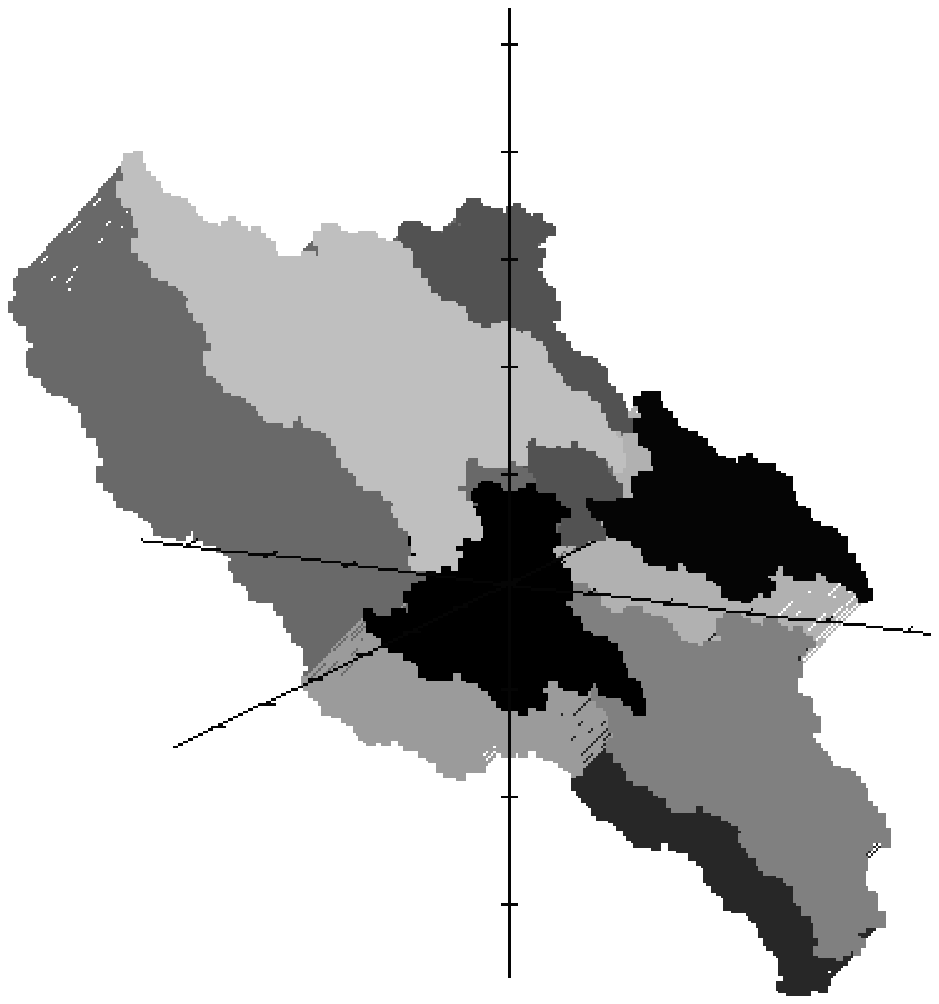}
\end{minipage}
\caption{The geometric representation  of  the  two-sided $\beta$-shift  for the Fi\-bo\-nacci-shift, the 
Tri\-bo\-nacci-shift,  and
the smallest Pisot number-shift.}\label{fig_markov}
\end{center}
\end{figure}

The geometric representation  map
$\widetilde{\varphi_{\beta}}$  of the golden ratio-shift
maps to ${\mathbb R}^2$. The ones for the Tribonacci number
and the smallest Pisot number map to ${\mathbb R}^3=\GC \times \GR$.
They are shown (up to a change of coordinates) in Fig. \ref{fig_markov}.
The sets  $\widetilde{{\mathcal R}_{\beta}}$
are unions of products of the different pieces of  the Rauzy fractal 
${\mathcal R}_{\beta}$ by finite real
 intervals of different heights.  For instance, in the Tribonacci case,
since the different pieces  in ${\mathcal R}_{\beta}$
correspond to  the sequences 
$(w_i)_{i\geq 0}$
 such that either 
$w_0=0$, or $w_0w_1=10$, or $w_0w_1=11$, this gives 
different constraints on  the sequences $(u_i)_{i\geq 1}$ which produce
the real component; for instance,
the piece which corresponds to the sequences $(w_i)_{i\geq 0}$
 such that 
 $w_0w_1=11$  implies the following constraint on  the sequences $(u_i)_{i\geq 1}$:
$u_0$ has to be equal to $0$.

\begin{remark}\label{rem:pi}
The  definition of the  representation map 
$\widetilde{\varphi_{\beta}}$ is inspired by \cite{sie_ari}
in the substitutive case. We cannot apply here directly
the substitutive formalism to the substitution $\sigma_{\beta}$ since it is  generally not 
Pisot :  the dominant eigenvalue  of $\sigma_{\beta}$ is still a Pisot number but other eigenvalues
may occur, as in the smallest Pisot case. The main difference here  is that we do not take those extra eigenvalues 
 into acount in our definition of $\widetilde{\varphi_{\beta}}$.
\end{remark}

\section{Geometric properties of the Rauzy fractal of  the $\beta$-numeration}\label{sec:pro}

The aim of this section is to prove some geometric properties of the Rauzy fractal
associated with  the $\beta$-numeration:
\begin{theorem}\label{theo:prop}
Let $\beta$ be a Pisot number.
\begin{enumerate}
\item The Rauzy fractal ${\mathcal  R}_{\beta}$ of the one-sided $\beta$-shift
has a  graph directed self-affine
structure in the sense of  \cite{mauldin}. More precisely:
\begin{itemize}
\item  it has a non-zero measure for
the Haar measure;
\item   there are $d$ pieces which  form the  self-affine structure,
where $d$ is equal to the  period  of  $d_{\beta}(1)^*$ if $\beta$
is a simple Parry number or to the sum of its preperiod   plus its period, otherwise.
\end{itemize}
\item  The Rauzy fractal $\widetilde{{\mathcal  R}_{\beta}}$ of the two-sided $\beta$-shift
has non-zero measure
for the Haar measure   $\mu_{\GK_{\beta}}$ 
over ${\mathbb K}_{\beta}\times \GR$. It is the disjoint union of $d$ cylinders obtained
as the product of each piece of the one-sided Rauzy fractal by a finite interval of $\GR$. 
\end{enumerate}
\end{theorem}
\subsection{Proof of Theorem \ref{theo:prop}}

The proof of the first point of this theorem is divided  into several steps; roughly, we  use
the self-affine structure of ${\mathcal  R}_{\beta}$  to deduce that
it has non-zero measure.
Similar statements are proved in the framework
 of Pisot  substitutions  in \cite{sie_ari}, but
we cannot use here directly  these statements  since the $\beta$-substitutions (Section \ref{sec:beta}) 
used here are not necessarily Pisot according to Remark \ref{rem:pi},
as for instance in the smallest
Pisot case.

\subsection*{Division of ${\mathcal  R}_{\beta}$ into subpieces.}
Let us first divide  (as illustrated in Figure \ref{fig_fractal})
the set ${\mathcal  R}_{\beta}$
into $d$ pieces where $d$ denotes the number of states in the minimal automaton
${\mathcal M}_{\beta}$ which recognizes 
the set of factors $F(X_{\beta})$ of the two-sided 
shift $X_{\beta}$.  

Let  $\widetilde{{\mathcal M}_{\beta}}$
denote the (non-deterministic) automaton obatined from 
 ${\mathcal M}_{\beta}$ by reversing the orientation of the edges. 
The set of sequences in  the left one-sided shift $X_{\beta}  ^l$
is  equal to the set of labels of infinite one-sided paths in the automaton $\widetilde{{\mathcal M}_{\beta}}$.
We label the states of the automaton  $\widetilde{{\mathcal M}_{\beta}}$  by $a_1,\cdots,a_d$. For
$i=1,\cdots,d$, one defines
 $${\mathcal  R}_{\beta}(i):=\left\{\begin{array}{ll}
 &\varphi_{\beta}((w_k)_{k\geq 0});
\  (w_k)_{k\geq 0}
 \in {\mathcal A}_{\beta}^{\GN^*}; 
(w_k)_{k \geq 0}  
\\ & \mbox{ is a path from the state }a_i \mbox{ in }\widetilde{ {\mathcal M}_{\beta}}
\end{array}\right\}.$$

Let us recall that $d$ is larger than or equal to the degree $r$ of the minimal polynomial
of $\beta$.  For instance,  $r=3$ and $d=5$ in the smallest Pisot case.
In particular, $d$ can be arbitrarily large when  $\beta$ is cubic, according to \cite{bassi}.

\subsection*{Self-affine decomposition.}
 We recall that  the map $h_{\beta}:\GK_{\beta} \rightarrow \GK_{\beta}$
stands for  the multiplication  map in ${\mathbb K}_{\beta}$
by the diagonal matrix whose diagonal  coefficients are given  by  $\delta_{\beta}(\beta)$ (Definition \ref{defi:h}).
Hence $\delta_{\beta}(\beta x)=h_{\beta}\delta_{\beta}(x)$ for every $x\in {\mathbb Q}(\beta)$. 

 Let us first  prove that for $i=1,\cdots,d$:
\begin{equation}\label{eq:ss}
{\mathcal  R}_{\beta}(i)=\cup _{1\leq j \leq d}\cup_{p, \   \sigma_{\beta}(j)=pis}
h_{\beta}({\mathcal  R}_{\beta}(j)) +\delta_{\beta}(|p|).
\end{equation}
 Let $i\in \{1,\cdots,d\}$ be given.
Let $(w_k)_{k\geq 0} \in {\mathcal  R}_{\beta}(i)$.
One has $$\begin{array}{ll}\varphi_{\beta}((w_k)_{k\geq 0})&=
\delta_{\beta}(\sum _{k\geq 1} w_k \beta ^k)+\delta_{\beta}(w_0)=
h_{\beta}\delta_{\beta}(\sum _{k\geq 1 } w_k  \beta ^{k-1})+\delta_{\beta}(w_0)\\
&=h_{\beta}\varphi_{\beta}((w_k)_{k\geq 1})
+\delta_{\beta}(w_0).
\end{array}$$

  Let $a_j$ denote the state in $\widetilde{{\mathcal M}_{\beta}}$
obtained by reading the label $w_0$ from the state $a_i$. By the definition of $\sigma_{\beta}$,
there exist $p=1^{w_0}$ and $s$ such that 
$\sigma_{\beta}(j)=pis$, that is, $w_0=|p|$.
Hence $$\varphi_{\beta}((w_k)_{k\geq 0})\in h_{\beta}({\mathcal  R}_{\beta}(j)) +\delta_{\beta}(|p|),
$$
which provides one inclusion  for the equality  (\ref{eq:ss}). The other 
inclusion is then immediate.

\subsection*{The Rauzy fractal has  non-zero measure}
The proof is an adaption of the corresponding proof
of \cite{sie_ari} which is done in the framework
of Pisot substitution dynamical systems; two properties 
are recalled below without a proof, the first one describing the action of the multiplication map
$h_{\beta}$  on  the Haar measure  $\mu_{\GK_{\beta}}$ of $\GK_{\beta}$:
\begin{lemma}[\cite{sie_ari}] \label{lem:mesu}
\begin{enumerate} 
\item For every Borelian set  $B$ of $\GK_{\beta}$ 
$$\mu_{\GK_{\beta}}(h_{\beta}(B))=\frac{1}{\beta} \mu_{\GK_{\beta}}(B).$$
\item Let ${\mathcal S}$  be a finite set included in $\GQ(\beta)$.
The set of points $\{\delta(P(1/\beta)); \ P \in 
{\mathcal S}[X]\}$ is a discrete set  in $\GK_{\beta}$.
\end{enumerate}
\end{lemma}

Let $(U_N)_{N\in\GN}$ denote  the linear  canonical  numeration system associated with $\beta$
according to \cite{BM89}, that is,
$U_0=1$,  and for all $k$, $U_k=t_1U _{k-1}+\cdots+t_kU_0+1$,
where $d_{\beta}^*(1)=(t_k)_{k\in \GN}$. 
Let us expand every integer  $i=0,1, \cdots, U_N-1$ in this system 
according to the greedy algorithm (this expansion being unique):
$$i=\sum_{0 \leq k \leq N-1} w_k ^{(i)}U_k.$$
According to \cite{BM89}, see also \cite{Loth},  the finite words    $w_{N-1}^{(i)}\cdots w_0^{(i)}$,
for $i=0,1, \cdots, U_N-1$  are all distinct and  all  belong to the set $F(X_{\beta})$ of factors of
elements of $X_{\beta}$.
Hence the sequences  $((w_0^{(i)}\cdots w_{N-1}^{(i)})0^{\infty} )$,
 for $i=0,1, \cdots, U_N-1$,
all belong to the left one-sided $\beta$-shift 
$X ^l _{\beta}$. 

Let ${\mathcal E}_N$ denote the image 
under the action of  $\varphi_{\beta}$ of this set of points.
One has for  $i=0,1, \cdots, U_N-1$, $$\varphi_{\beta}((w_0^{(i)}\cdots w_{N-1}^{(i)})0^{\infty} )
=\delta_{\beta}(\sum_{0 \leq k \leq N-1} w_k ^{(i)} \beta ^i).$$
The points in  ${\mathcal E}_N$
are all distinct since  for $i=0,1, \cdots, U_N-1$, the finite words    $w_{N-1}^{(i)}\cdots w_0^{(i)}$
  are all distinct. There thus exists $B >0$ such that
$\mbox{Card } {\mathcal E}_N  =U_N> B \beta ^N$, since $\beta$ is a Pisot number.

Let  us apply  now Lemma \ref{lem:mesu} with  
${\mathcal S}=\{0,1,\cdots,[\beta]\}$. 
Hence there exists a constant  $A>0$
such that the distance between two elements
of ${\mathcal E}_N$ is larger than $A$.
Let us define now for every non-negative integer  $N$, the set $${\mathcal B}_N=\cup_{z \in {\mathcal E}_N}
(h_{\beta}) ^N B(z,A/2),$$ where 
$B(z,A/2)$  denotes the closed ball in
$\GK_{\beta}$ of center $z$ and radius $A/2$.
According to Lemma \ref{lem:mesu}
and to the fact  that 
$\mbox{Card } {\mathcal E}_N > B \beta ^N$, for all $N$, there exists
$C>0$ such that $\mu_{\GK_{\beta}}({\mathcal B}_N)>C$, for all $N$.

The main point now is that the sequence of compact sets  $({\mathcal B}_N)_{N \in \GN}$
 converges toward  a subset of ${\mathcal R}_{\beta}$  with respect  to the 
Hausdorff metric; indeed   for a fixed positive integer $p$, for $N$ large enough,
${\mathcal B}_N \subset {\mathcal R}_{\beta}(1/p):=
\{x \in \GK_{\beta}; \ d(x,{\mathcal R}_{\beta}) \leq 1/p\}$. 
Independently, fix $\varepsilon >0$; since the sequence $({\mathcal R}_{\beta}(1/p))_p$ converges toward  ${\mathcal R}_{\beta}$,
 there exists $p >0$ such that
$\mu_{\GK_{\beta}}({\mathcal R}_{\beta}(1/p)) \leq \mu_{\GK_{\beta}}({\mathcal R}_{\beta})
+\varepsilon$. This finally 
 implies $\liminf \mu_{\GK_{\beta}}({\mathcal B}_N) \leq  \mu_{\GK_{\beta}}({\mathcal R}_{\beta})
+\varepsilon.$ Since this holds for every $\varepsilon>0$, one obtains $
 \mu_{\GK_{\beta}}({\mathcal R}_{\beta})\geq C > 0$.

\subsection*{Computation of the measures and self-affine structure}
Let us  prove now that the union
in (\ref{eq:ss}) is a disjoint union up to sets of zero measure.
One has for a given  $i\in\{1,\cdots,d\}$ according to (\ref{eq:ss}) and to Lemma
\ref{lem:mesu}
\begin{equation}\label{eq:mes3}
\begin{array}{ll}
\mu_{\GK_{\beta}}({\mathcal  R}_{\beta}(i)) &\leq \sum_{j=1,\cdots,d, \  \sigma_{\beta}(j)=pis}\mu_{\GK_{\beta}}(h_{\beta}({\mathcal  R}_{\beta}(j)) )\\
&\leq 1/\beta\sum_{j=1,\cdots,d, \  \sigma_{\beta}(j)=pis}\mu_{\GK_{\beta}}({\mathcal  R}_{\beta}(j) ).
\end{array}
\end{equation}
Let ${\bf m}=(\mu_{\GK_{\beta}}({\mathcal  R}_{\beta}(i)))_{i=1,\cdots, d}$ denote the vector in $\GR^d$
of measures in $\GK_{\beta}$ of  the pieces of the Rauzy fractal. We know  from what precedes that ${\bf m}$
is a   non-zero vector with non-negative entries.
According to Perron-Frobenius theorem, 
the previous equality implies that ${\bf m}$
is  an eigenvector
  of  the primitive incidence  matrix of  the substitution ${\sigma}_{\beta}$. We thus have equality in (\ref{eq:mes3})
which implies that the unions are disjoint up to sets of zero measure.
One similarly proves that this equality in measure still holds by replacing
$\sigma$ by $\sigma ^n$,
for every $n$.

Now take two distinct pieces, $({\mathcal  R}_{\beta}(j))$ and $({\mathcal  R}_{\beta}(k))$ say,
with $j\neq k$. There exists $n$ such that both $\sigma^n (j)$ and $\sigma ^n (k)$ admit as first letter $1$. Hence they both occur in (\ref{eq:ss}) for $i=1$ with the same translation term
(which is indeed equal to $0$)  and they are thus distinct. We hence have proved
that the $d$ pieces of the Rauzy  fractal $({\mathcal  R}_{\beta}(j))$
are disjoint up to sets of zero measure, which ends  the proof of Assertion (i), that is, ${\mathcal  R}_{\beta}$ has
a self-affine structure.

\subsection*{$\widetilde{{\mathcal  R}_{\beta}}$
has non-zero measure}
It remains now to prove the second point of the theorem.
The assertion that $\widetilde{{\mathcal  R}_{\beta}}$
has non-zero measure is  a direct consequence of the structure of the Rauzy fractal $\widetilde{{\mathcal  R}_{\beta}}$
studied above since $\widetilde{{\mathcal  R}_{\beta}}$ can be decomposed  as  the 
disjoint (in measure) union
of $\widetilde{{\mathcal  R}_{\beta}}(i)$, $i=1,\cdots,d$, where
 $$\widetilde{{\mathcal  R}_{\beta}}(i):=\left\{\begin{array}{ll}
 &\widetilde{\varphi_{\beta}}(u,w);
\ (u,w)  \mbox{ is a two-sided path }\mbox{ in }\widetilde{ {\mathcal M}_{\beta}}\\
&\mbox{ such that }w \mbox{ starts from the state }a_i
\end{array}\right\},$$ which the product 
of ${\mathcal  R}_{\beta}(i)$ by a finite  interval of $\GR$ of non-zero measure.

\subsection{Examples} \label{sec:ex2}
Let us pursue the  study of  three of the examples of Section  \ref{sec:ex}.

\subsubsection*{The Tribonacci number} 
Let us recall that $${\mathcal R}_{\beta}=\{\sum_{i\geq 0}
w_i \alpha^i; \  \forall i, \ w_i \in \{0,1\}, \ w_i w_{i+1}w_{i+2}\neq 0\}.$$
One easily checks thanks to the automaton $\widetilde{{\mathcal M}_{\beta}}$  shown in Fig. \ref{exemplesAuto} that ${\mathcal R}_{\beta} (1)$  corresponds to the sequences  $(w_i)_{i\geq 0}$ such that
$w_0=1$, 
${\mathcal R}_{\beta} (2)$  corresponds to the  set of sequences  $(w_i)_{i\geq 0}$
 such that $w_0w_1=10$, and lastly, ${\mathcal R}_{\beta} (3)$ corresponds to the  sequences  $(w_i)_{i\geq 0}$
 such that $w_0w_1=11.$
One has 
$$\left\{
\begin{array}{l}
{\mathcal R}_{\beta} (1)=\alpha({\mathcal R}_{\beta} (1)\cup {\mathcal R}_{\beta} (2)\cup{\mathcal R}_{\beta} (3))\\
{\mathcal R}_{\beta} (2)=\alpha({\mathcal R}_{\beta} (1))+1\\
{\mathcal R}_{\beta} (3)=\alpha({\mathcal R}_{\beta} (2))+1.
\end{array}
\right.
$$

\subsubsection*{The smallest Pisot  number} 
One has $$\begin{array}{ll}
{\mathcal R}_{\beta}=&
\{
\sum_{i\geq 0}
w_i \alpha^i; \  \forall i, \ w_i \in \{0,1\}, \  \mbox{ if } w_i =1,\\
& \mbox{ then }w_{i+1}=w_{i+2}=
w_{i+3}=w_{i+4}= 0\}.
\end{array}$$
One easily checks that ${\mathcal R}_{\beta} (1)$  corresponds to the sequences  $(w_i)_{i\geq 0}$ such that
$w_0=0000$, 
${\mathcal R}_{\beta} (2)$  corresponds to the  set of sequences  $(w_i)_{i\geq 0}$
 such that $w_0w_1=1$,  ${\mathcal R}_{\beta} (3)$ corresponds to the  sequences  $(w_i)_{i\geq 0}$
 such that $w_0w_1=01,$ ${\mathcal R}_{\beta} (4)$ corresponds to the  sequences  $(w_i)_{i\geq 0}$
 such that $w_0w_1=001,$ 	and lastly 
${\mathcal R}_{\beta} (5)$ corresponds to the  sequences  $(w_i)_{i\geq 0}$
 such that $w_0w_1=0001.$

One has $$\left\{
\begin{array}{l}
{\mathcal R}_{\beta} (1)=\alpha({\mathcal R}_{\beta} (1)\cup {\mathcal R}_{\beta} (5))\\
{\mathcal R}_{\beta} (2)=\alpha({\mathcal R}_{\beta} (1))+1\\
{\mathcal R}_{\beta} (3)=\alpha({\mathcal R}_{\beta} (2))\\
{\mathcal R}_{\beta} (4)=\alpha({\mathcal R}_{\beta} (3))\\
{\mathcal R}_{\beta} (5)=\alpha({\mathcal R}_{\beta} (4)).
\end{array}
\right.
$$
\subsubsection*{The $(2+\sqrt 2)$-shift} 
\begin{figure}[ht]
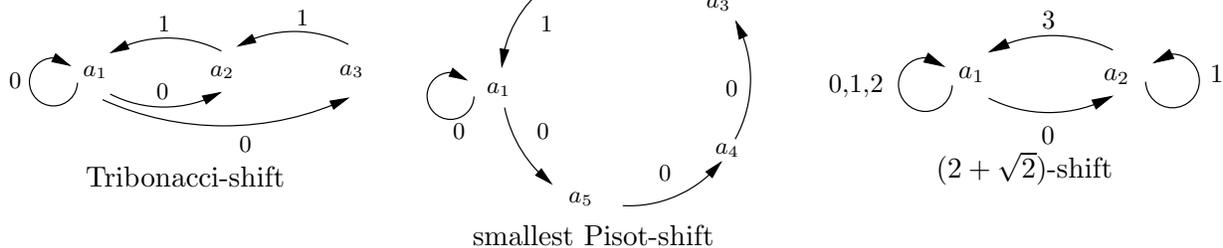

\begin{center}
\begin{minipage}{5cm}
\begin{center}
\input{auto_Tribo.pstex_t}

Tribonacci-shift
\end{center}
\end{minipage}
\begin{minipage}{5.6cm}
\begin{center}
\input{auto_PPP.pstex_t}

smallest Pisot-shift
\end{center}\end{minipage}
\begin{minipage}{5.6cm}
\begin{center}
\input{auto_NU.pstex_t}

$(2+\sqrt 2)$-shift
\end{center}\end{minipage}
\caption{Reversed minimal automaton $\widetilde{{\mathcal M}_{\beta}}$ describing the structure
of the $\beta$-shift.}\label{exemplesAuto}\end{center}
\end{figure}

One has $${\mathcal R}_{\beta}=
\{
\sum_{i\geq 0}
w_i (2-\sqrt 2)^i; \  \forall i, \ w_i \in \{0,1,2\}, \ (w_i)_{i\geq 0} \leq_{\mbox{lex}} (31^{\infty})\}.
$$
In this non-simple Parry case, we cannot  express as easily as previously
the sets ${\mathcal R}_{\beta} (1)$
and ${\mathcal R}_{\beta} (2)$: one checks in Figure \ref{exemplesAuto} that 
there exist cycles from $a_1$ to $a_1$
and from  $a_2$ to $a_2$, which  implies that both  ${\mathcal R}_{\beta} (1)$ and
${\mathcal R}_{\beta} (2)$ contain sequences with 
arbitrarily long common prefixes, such as $1^n$ for every $n$. 

One has $$\left\{
\begin{array}{ll}
{\mathcal R}_{\beta} (1)=& (2-\sqrt 2)({\mathcal R}_{\beta} (1))
\cup ((2-\sqrt 2) ({\mathcal R}_{\beta} (1))+1)
\cup ((2-\sqrt 2)({\mathcal R}_{\beta} (1))+2)\\
&\cup (2-\sqrt 2)({\mathcal R}_{\beta} (2))\\
{\mathcal R}_{\beta} (2)=& (2-\sqrt 2)({\mathcal R}_{\beta} (1))+3+ (2-\sqrt 2)({\mathcal R}_{\beta} (2))+1.\end{array}
\right.
$$

\section{Characterization of purely periodic points}\label{sec:car}
We can now state the main theorem of this paper.

\begin{theorem}\label{theo:main2}
Let $\beta$ be a Pisot number.
For all $x \in {\mathbb Q}(\beta) \cap [0,1)$,
the $\beta$-expansion of  $x$ is purely periodic
if and only if $(\delta_{\beta}(x), x) \in
\widetilde{\mathcal R_{\beta}}=\widetilde{\varphi_{\beta}} ({X_{\beta}})$.
\end{theorem}

\prv

Let us assume that the $\beta$-expansion of  $x$ 
is purely periodic. 
Write $x$ as $x= 0.(a_1 \dots a_L)^{\infty}$,
and set
$w= \!\!^{\infty} (a_1 \dots a_L)$ and $u= (a_1 \dots a_L)^{\infty}$
(that is, $w=(w_i)_{i\geq 0}$ with $w_0\cdots w_{L-1}=a_L \dots a_1$,
and $w_{i+L}=w_i$ for all $i$,  and $u=(u_i)_{i\geq 1}$, with $u_1\cdots u_{L}=a_1 \dots a_L$,
and $u_{i+L}=u_i$ for all $i$).
One has $(w,u)  \in {X_{\beta}}$ acording to (\ref{eq:real}).
Let us compute $\widetilde{\varphi_{\beta}}(w,u)$.  Note
that the second  coordinate of $\widetilde{\varphi_{\beta}}(w,u)$ is
$x$ by definition:
$$\sum_{i\geq 1} u_i\beta^{-i} =
\frac{a_1 \beta^{-1}  + \dots + a_L \beta ^{-L}}{1-\beta^{-L}} = 
 \frac{a_1 \beta^{L-1} + \dots + a_L}{\beta^L-1} = x.$$

 Futhermore,
 $\lim_{n \to \infty} \delta_{\beta}(\beta^n) = 0$
in $\GK_{\beta}$. We thus have
\begin{eqnarray*}
\widetilde{\varphi_{\beta}}(w,u)&=&\hspace{-0.1cm}  \left( 
-\delta_{\beta}(\sum_{i\geq 0} w_i\beta^{i}),\sum_{i\geq 1} u_i\beta^{-i } \right)\\
&=&
  \left(- \lim_{n\to \infty} \delta_{\beta}
((a_L+\dots+a_1 \beta^{L-1} )(1+\beta^L+\dots + \beta^{nL})),x\right) \\
&=&
 \left( \lim_{n\to \infty} \delta_{\beta}
\left(-(a_1 \beta^{L-1} + \dots + a_L)\frac{1-\beta^{nL}}{1-\beta^L}\right),x\right)\\
&=&
 \left(\delta_{\beta}
\left(\frac{a_1 \beta^{L-1} + \dots + a_L}{\beta^L-1}\right),x\right)
\\ &=& (\delta_{\beta}(x),x),
\end{eqnarray*}
hence $(\delta_{\beta}(x),x)\in \widetilde{{\mathcal R}_{\beta}}$.

\medskip

Consider now the reciprocal and let  $x\in {\mathbb Q}(\beta)$ 
such that   $(\delta_{\beta}(x),x) \in \widetilde{{\mathcal R}_{\beta}}$.

The sketch of the proof is the following and is inspired by
a similar discrete argument in
\cite{Hui2}, dealing with the case $\beta$ is a unit.
We  introduce  below a finite subset ${\mathcal T}_x$
of $\widetilde{{\mathcal R}_{\beta}}$, which depends on
$x$, 
and which is  stable
under the action of the map
$\widetilde{T}_{\beta}$ which is onto on it.
In order to define this finite set, we 
 take into account 
all the  ${\mathcal I}$-adic topologies
which correspond to  prime ideals ${\mathcal I}$
which do not appear in  the decomposition
(\ref{decompo}).
These topologies are the extensions
on ${\mathbb Q}(\beta)$  of $p$-adic topologies for
 the primes $p$ which do not divide
the constant term
of the minimal polynomial $P_{\beta}$ of $\beta$ according to  \cite{sie_ari}.
Roughly speaking, one introduces further restrictions
over $\GQ(\beta)$
which involve the primes that
were not already considered
in $\GK_{\beta}$.
We will  prove  below that the set ${\mathcal T}_x \subset \widetilde{{\mathcal R}_{\beta}}
\subset
\GK_{\beta} \times \GQ(\beta)$
 is finite  by using  the fact that the second  coordinate
of its elements, which belongs
to $\GQ(\beta)$,   is bounded  for all the topologies on
$\GQ(\beta)$.
\medskip

Let us first introduce  the following set ${\mathcal S}_x$:
$$ {\mathcal S}_x = \{z  \in {\mathbb Q}(\beta), \
|z|_{\mathcal I}\leq \mbox{max}(\,|x|_{\mathcal I},1 \,), \,
 \mbox{ for every prime ideal } \,
 {\mathcal I} \, \mbox{ such that} \, |\beta|_{\mathcal I}=1\}.
$$

One has
\begin{itemize}
\item $x \in {\mathcal S}_x$ since
 $|x|_{\mathcal I}\leq \mbox{max}(\,|x|_{\mathcal I},1 \,)$
for every prime ideal ${\mathcal I}$. 
\item ${\mathbb Z}\subset {\mathcal S}_x$ since
for every integer $N$ and for every ideal ${\mathcal I}$,
one has $|N|_{\mathcal I} \leq 1 \leq \mbox{max}(\,|x|_{\mathcal I},1 \,)$.
\item $\beta {\mathcal S}_x \subset {\mathcal S}_x$ and
 $\beta^{-1} {\mathcal S}_x\subset {\mathcal S}_x$ since
$|\beta z|_{\mathcal I}= |\beta|_{\mathcal I}\, | z|_{\mathcal I}
= | z|_{\mathcal I}$  if $|\beta|_{\mathcal I}=1$,
and as well $|\beta^{-1} z|_{\mathcal I}= | z|_{\mathcal I}$.
\item $({\mathcal S}_x,+)$ is a group since 
the ${\mathcal I}$-adic valuations are ultrametric.
\end{itemize}

Hence one deduces that:

\begin{itemize}
\item ${\mathcal S}_x$ is stable under the action of  $T_{\beta}$; indeed 
$$T_{\beta} {\mathcal S}_x \subset \beta {\mathcal S}_x
- \GZ \subset {\mathcal S}_x - {\mathbb Z}
\subset {\mathcal S}_x - {\mathcal S}_x \subset {\mathcal S}_x.$$

\item For every integer $N$,
 one has $\beta^{-1}({\mathcal S}_x + N) \subset \beta^{-1}{\mathcal S}_x \subset {\mathcal S}_x$. 
\end{itemize}

In the original proof of
\cite{Hui2}, the analogous of the set
${\mathcal S}_x$
is $\GZ[\beta]/q$
for an integer $q$ which depends on $x$, 
this set being no more stable under the  multiplication
by $1/\beta$ in the non-unit case.
The keypoint of the present proof  is that ${\mathcal S}_x$ is not only stable under the action
of $T_{\beta}$ but also under
the multiplication by $1/\beta$, even when  $\beta$ is not a unit.

\medskip 

Let us consider now the following subset of $\widetilde{{\mathcal R}_{\beta}}$
obtained by first embeding the points
of ${\mathcal S}_x$ into
${\mathbb K}_{\beta}\times \GR $, and then intersecting it
with  the bounded set $\widetilde{{\mathcal R}_{\beta}}$:
 $${\mathcal T}_x=(\delta_{\beta},Id)({\mathcal S}_x) \cap \widetilde{{\mathcal R}_{\beta}}\subset \GK_{\beta} \times \GQ(\beta).$$
Let us first observe that $(\delta_{\beta}(x),x)\in {\mathcal T}_x$.

\medskip

\begin{itemize}
\item  {\bf The set ${\mathcal T}_x$ is finite}. 
Indeed let us prove that
all the absolute values of $\GQ(\beta)$ take  bounded  values
on the last  coordinate
of the elements of ${\mathcal T}_x$.
Let $z\in {\mathcal S}_x$ such that
$(\delta_{\beta}(z),z)\in \widetilde{{\mathcal R}_{\beta}}$.
Note first that the usual absolute value of $z$ is bounded.
Furthermore
\begin{itemize}
\item if  $|\cdot |$ is
an Archimedean valuation or if
$|\cdot |_{\mathcal I}$ 
is an ultrametric valuation satisfying $|\beta|_{\mathcal I}\not=1$,
this valuation appears by construction
in ${\mathbb K}_{\beta}$,  and
 thus $|z|$ is bounded
since $(\delta_{\beta}(x),x)$ belongs to the bounded set
$\widetilde{{\mathcal R}_{\beta}}$.
\item if
$|\cdot |_{\mathcal I}$ 
is an ultrametric valuation satisfying
$|\beta|_{\mathcal I}=1$, then
$|z|_{\mathcal I}$ is bounded 
 by definition of ${\mathcal S}_x$. 
\end{itemize}

\item {\bf  The map
$\widetilde{T_{\beta}}$ is onto over ${\mathcal T}_x$.} Let 
 $(\delta_{\beta}(z),z) \in  {\mathcal T}_x$, with $z \in {\mathcal S}_x$.
 By definition,
there exists $(w,u)=((w_i)_{i\geq 0},(u_i)_{i\geq 1}) \in \widetilde{X_{\beta}}$ such  that
 $\widetilde{\varphi_{\beta}}(w,u)=(\delta_{\beta}(z),z)$. 
Hence $z=\sum_{i\geq 1} u_i \beta^{-i}$ 
and $\delta_{\beta}(z)=-\delta_{\beta}(\sum_{i\geq 0} w_i \beta^{i})$.
An easy computation shows that
\begin{equation} \label{eq:minus}
\widetilde{\varphi_{\beta}} \circ S^{-1}(w,u)=
\left(\delta_{\beta}\left(\frac{z+w_0}{\beta}\right), 
\frac{z+w_0}{\beta}\right).
\end{equation} 
Let us note that this computation works thanks to the introduction of a sign minus
before $\delta_{\beta}$
in the expression of $\widetilde{\varphi_{\beta}}$.

Let us first assume that $w_0u_1u_2 \dots\not=d^*_\beta(1)$. 
One thus deduces  according to(\ref{eq:prop1}),   that
\begin{align*}
 (\delta_{\beta}(z),z) = \widetilde{\varphi_{\beta}}(w,u)
= \widetilde{\varphi_{\beta}} \circ S \circ S^{-1}(w,u)
&= \widetilde{T_{\beta}} \circ \widetilde{\varphi_{\beta}} \circ S^{-1}(w,u)
 &= \widetilde{T_{\beta}} \circ (\delta_{\beta},Id) \left(\frac{z+w_0}{\beta}\right).
\end{align*}
On the one hand $(\delta_{\beta},Id) \left(\frac{z+w_0}{\beta}\right)
= \widetilde{\varphi_{\beta}}\circ S^{-1}(w,u) \in \widetilde{{\mathcal R}_{\beta}}$, and on the other hand
$\frac{z+w_0}{\beta} \in \frac{{\mathcal S}_x+w_0}{\beta}
\subset \beta^{-1}({\mathcal S}_x+{\mathbb Z}) \subset{\mathcal S}_x$.  Hence,
$ (\delta_{\beta}(z),z) \in \widetilde{T_{\beta}}( {\mathcal T}_x).$

\smallskip

Assume now that $w_0u_1u_2 \dots\not=d^*_\beta(1)$. By definition $w_0 = [\beta]$.
Moreover,  one checks that $z=\beta -w_0$. One has
$$ \widetilde{T_\beta}\circ(\delta_\beta,Id)(1) = (h_\beta(1)-[\beta]\delta_{\beta}(1),\beta-[\beta])
=(\delta_{\beta}(\beta)-\delta_{\beta}([\beta]), \beta-[\beta])
=(\delta_\beta(z),z).$$ 
Since  $\frac{z+w_0}{\beta}$ and from (\ref{eq:minus}), one has $(\delta_\beta,Id)(1) \in 
\widetilde{{\mathcal R}_{\beta}}$ so that $ (\delta_{\beta}(z),z) \in \widetilde{T_{\beta}}( {\mathcal T}_x).$

\item {\bf The set
${\mathcal T}_x$ is stable by
 $\widetilde{T_{\beta}}$.}
Let 
 $(\delta_{\beta}(z),z) \in  {\mathcal T}_x$, with $z \in {\mathcal S}_x$.
 By definition,
there exists $(w,u)=((w_i)_{i\geq 0},(u_i)_{i\geq 1}) \in \widetilde{X_{\beta}}$ such  that
 $\widetilde{\varphi_{\beta}}(w,u)=(\delta_{\beta}(z),z)$. 
One has  according to (\ref{eq:prop2}) 
$$\widetilde{T_{\beta}}(\delta_{\beta}(z),z) = (\delta_{\beta},Id)\circ T_{\beta}( z).$$

Let us first suppose that $u\not=d^*_\beta(1)$. One deduces from 
(\ref{eq:prop1})
$$\widetilde{T_{\beta}}\circ \widetilde{\varphi_{\beta}}(w,u)=
 \widetilde{\varphi_{\beta}}\circ S ( w,u ).$$

Now if $u=d^*_\beta(1)$, then $z= \sum_{i\geq 1}u_i \beta^{-i}= 1$, $T_\beta(z)=\beta-u_1$ and $ \delta_{\beta}(\sum_{i\geq 0}w_i \beta^{i}=\delta_\beta(1)$ so that
$$ \widetilde{\varphi_{\beta}}\circ S ( w,u )=(-\delta_{\beta}(u_1+\sum_{i\geq 0}w_i \beta^{i+1}),
\sum_{i\geq 2} u_i \beta^{-i+1}) 
= (-\delta_{\beta}(u_1)+\delta_{\beta}(\beta), \beta-u_1)=(\delta_{\beta},Id)\circ T_{\beta}( z).$$
\end{itemize}

\medskip

Hence, since $\widetilde{T_{\beta}}$ is onto over a stable  finite set, it is 
one-to-one over this set
${\mathcal T}_x$.
By construction, 
   $(\delta_{\beta}(x),x)
\in {\mathcal T}_x$. Hence there exists 
an integer $n$ such that
$$(\delta_{\beta}(x),x) = \widetilde{T_{\beta}}^n(\delta_{\beta}(x),x)
= (\delta_{\beta}({T_{\beta}}^n(x)), {T_{\beta}}^n(x)).$$
We thus deduce that the expansion of   $x$
is purely periodic.
\fprv

\section{Conclusion}

One can also introduce a notion of generalized Rauzy fractal in the non-Pisot case.
Namely, a  similar characterization holds in the non-Pisot case by introducing
in the representation map
$\widetilde{\varphi_{\beta}}$
the
coordinate
$\sum_{i\geq 1} u_i \lambda ^{-i}$
for the conjugates $\lambda$ of modulus strictly larger that $1$ (see for instance
\cite{kenyon} for a similar description).
Yet we are not able to prove that the union in  Equation (\ref{eq:ss})
is disjoint up to sets of zero measure and that the Haar measure of the generalized Rauzy
fractal  is non-zero.

It remains now to use this formalism to study  further  topological 
or metrical properties 
of the  sets ${\mathcal R}_{\beta}$
 and $\widetilde{{\mathcal R}_{\beta}} $, even in the non-Pisot case, which will be the 
object of a subsequent paper.
Our motivations are the following:   first   the construction of  explicit Markov partitions of 
endomorphisms of the torus as initiated in 
\cite{these_anne},   second,  the study of rational numbers having a   purely
periodic expansion  in the flavour of \cite{AKI1,schmidt80}, and third,
the spectral study of $\beta$-shifts in the Pisot non-unit case according to \cite{sie_ari}.

\bibliography{shift_court}
\bibliographystyle{alpha}

\end{document}